\documentclass[11pt]{amsart}
\usepackage[]{amsmath,amssymb,fullpage}

\begin{document}
\title [QD and of crossed products]  {Quasidiagonality of crossed products}
\author [S. Orfanos]  {Stefanos Orfanos}
\address  {Department of Mathematical Sciences, University of Cincinnati, Cincinnati, OH, 45221}
\email {stefanos.orfanos@uc.edu}

\begin{abstract}
We prove that the crossed product $A\rtimes_\alpha G$ of a separable unital quasidiagonal C*-algebra $A$ by a discrete countable amenable maximally almost periodic group $G$ is quasidiagonal, provided that the action $\alpha$ is almost periodic. This generalizes a result of Pimsner and Voiculescu.
\end{abstract}

\maketitle

\theoremstyle{plain}
\newtheorem {thm} {Theorem} 
\newtheorem {lem} [thm]{Lemma} 
\newtheorem {cor} [thm]{Corollary}
\newtheorem {pro} [thm]{Proposition}

\theoremstyle{definition}
\newtheorem {dfn} {Definition}
\newtheorem {que} {Question}
\newtheorem* {ex} {Example}

\theoremstyle{remark}
\newtheorem* {rem} {Remark}
\newtheorem* {akn} {Aknowledgments}
\newtheorem* {pr} {Proof}

\section*{Introduction}
Finite dimensional approximation properties of C*-algebras have been instrumental in unraveling the structure of these algebras and also in our efforts to classify them. Quasidiagonality is one of many such approximation properties. Informally speaking, it amounts to the algebra having local approximate embeddings into matrix algebras (Voiculescu's abstract characterization; see \cite{5} for the exact statement). The study of quasidiagonality was initiated by Halmos; he defined quasidiagonal operators as a generalization of the block-diagonal operators. Quasidiagonality has also topological implications; note that it is homotopy invariant (\cite{5}). For connections of quasidiagonality with other properties of C*-algebras and with the classification program, we refer the reader to \cite{2}.

In this note we prove that crossed products of unital, separable, quasidiagonal C*-algebras by discrete countable amenable maximally almost periodic groups are quasidiagonal, given that the group action is almost periodic. The subgroup structure of such groups allows us to construct certain tilings and from them, a sequence of projections that implement quasidiagonality. The crossed products are then treated in a similar manner as in the proof of M. Pimsner and D. Voiculescu, who showed in \cite{4}, among other things, that the crossed product of a unital separable quasidiagonal C*-algebra by the integers is quasidiagonal, if the action is almost periodic. 

An application of our result to a broad class of crossed products, the \emph{generalized Bunce--Deddens algebras}, is discussed in a subsequent paper (\cite{3}).

\section{Background and construction of the projections $P_n$}
We first recall some known definitions and facts. All groups hereafter will be discrete and countable and all C*-algebras unital and separable.

\begin{dfn}
 The group $G$ is \emph{residually finite} if it has a separating family of finite index normal subgroups.
\end{dfn}

\begin{lem}
$G$ is residually finite if and only if for any finite set $F$ of $G$, there is a finite index normal subgroup $L$ of $G$ such that $F^{-1}F \cap L=\{ e\}$.
\end{lem}

\begin{pr}
For all integers $0\le i<j\le m$, consider $H_{ij} \lhd G$ of finite index, with the property $x_i^{-1}x_jH_{ij}\neq H_{ij}$. The intersection $L = \bigcap _{i<j} H_{ij}$ is a finite index normal subgroup of $G$. We claim that $F^{-1}F \cap L=\{e\}$. Assume not, then there are $x_i, x_j$ in $F$ such that $x_i^{-1}x_j\in L\subset H_{ij}$ which forces $x_i^{-1}x_j = e$. The converse is clear.\qed 
\end{pr}

\begin{dfn} 
A group $G$ is \emph{amenable} if there exists a sequence $e\in F_1\subset F_2\subset\dotsc\subset F_n\dotsc$ of finite sets of $G$ such that:
\begin{enumerate}
\item
$\bigcup_{n\ge 1} F_n=G$
\item
$\lim_{n\to\infty}\frac{|F_n\triangle F_ns|}{|F_n|}=0$ for all $s\in G$
\end{enumerate}
The sequence $(F_n)_n$ is called a \emph{sequence of (right) F\o lner sets} for $G$.
\end{dfn}

A \emph{tiling} of $G$ is a covering of $G$ by disjoint (right) translates of a single set, the \emph{tile}. In the lemma below, $K_n$ is the tile and $L_n$ is the set of \emph{tiling centers} of a tiling of $G$, for every $n\ge 1$. B. Weiss has constructed tiles of discrete countable amenable residually finite groups $G$ which are also F\o lner sets (see \cite{6}). For our purposes, it suffices to know that they contain F\o lner sets.

\begin{lem}
Suppose $G$ is an amenable and residually finite group. Then for every $n=1,2,\dotsc$ the group $G$ has a tiling of the form $G=K_nL_n$, where for each $n$, $K_n$ is a finite subset of $G$ containing the F\o lner set $F_n$ and $L_n$ is a normal subgroup of finite index in $G$.
\end{lem}

\begin{pr}
For each $n=1,2,\dotsc$ consider the F\o lner set $F_n\subset G$. Use the previous lemma to get $L_n\lhd G$ such that $F_n^{-1}F_n \cap L_n=\{ e\}$. By adding at most finitely many elements to $F_n$, get a set $K_n\supset F_n$ such that $K_n^{-1}K_n \cap L_n=\{ e\}$ and moreover $|K_n|=\lbrack G:L_n\rbrack$. Then $G=K_nL_n$ \qed
\end{pr}

Recall that the \emph{left regular representation} of a group $G$ is a unitary representation $\lambda :G \to \mathcal{U}(\ell^2(G))$ given by $\lambda (s)\delta _x=\delta _{sx}$, $s\in G$. We obtain the \emph{reduced C*-algebra} $C_r^*(G)$ of the group $G$ by taking the closure of the span of $\{\lambda (s): s\in G\}$. The \emph{full C*-algebra} $C^*(G)$ is obtained as the closure of the universal representation. Moreover, recall that an \emph{action} of $G$ on a C*-algebra $A$ is a group homomorphism $\alpha$ from $G$ into the group of *-automorphisms of $A$. Then one considers the *-algebra $C_c(G,A)$ of finitely supported functions on $G$ with values in $A$, where multiplication is defined as the $\alpha$-twisted convolution. A \emph{covariant representation} is a triple $(\sigma ,u,\mathcal{H})$, where $\sigma :A\to \mathcal{B(H)}$ is a *-representation and $u:G\to \mathcal{B(H)}$ a unitary representation, such that $u(s)\sigma (a)u(s)^* = \sigma (\alpha (s)a)$ for every $s\in G$ and $a\in A$. An example: if $\rho :A\to \mathcal{B(H)}$ is a faithful *-representation and $\lambda$ is the left regular representation of $G$, define $\sigma :A\to \mathcal{B(H}\otimes\ell^2(G))$ to be $\sigma (a)(h\otimes \delta_x) =\rho (\alpha(x^{-1})a)h\otimes\delta_x$. Then the triple $(\sigma ,I_{\mathcal{H}}\otimes \lambda, \mathcal{H}\otimes\ell^2(G))$ is a covariant representation which we call the \emph{regular representation}. Covariant representations correspond precisely to *-representations of $C_c(G,A)$. The completion with respect to the norm coming from the regular representation gives the \emph{reduced crossed product} $A\rtimes _{\alpha ,r}G$, and the one coming from the universal representation gives the \emph{(full) crossed product} $A\rtimes _{\alpha}G$. The two norms coincide whenever $G$ is an amenable group. If in addition $A=\mathbb{C}$, we get $C_r^*(G)=C^*(G)$.

\begin{dfn}
A linear operator $T$ on a separable Hilbert space $\mathcal{H}$ is \emph{quasidiagonal} if there exists a sequence of finite rank self-adjoint orthogonal projections $P_n$ in $\mathcal{B(H)}$ satisfying:
\begin{enumerate}
\item $ P_n\to I_{\mathcal{H}}$ as $n\to\infty$, and
\item $\|[T,P_n]\| \to 0$ as $n\to\infty$ 
\end{enumerate}
A separable set of operators $\mathcal{A}$ is \emph{quasidiagonal} if every operator $T$ in a set of dense linear span in $\mathcal{A}$ is quasidiagonal with respect to the same sequence $(P_n)_n$. An abstract C*-algebra $A$ is \emph{quasidiagonal} if it has a faithful representation to a quasidiagonal set of operators.
\end{dfn}

We now turn to the construction of the projections $P_n$. First, define the function $\phi _n :G\to\left[ 0,1\right]$ by \[ \phi _n(x)=\sqrt{\frac{|K_n\cap F_nx|}{|F_n|}}\mbox{, }x\in G\] One observes that supp$\phi _n=F_n^{-1}K_n$ is finite for all $n\ge 1$.

\begin{lem}
The function $\phi _n$ enjoys the following properties:
\[\sum_{x\in\ yL_n} \phi _n^2(x)=1 \mbox{ and }\sum_{x\in yL_n} |\phi _n^2(x)-\phi _n^2(sx)|\leq \frac{|F_n\triangle F_ns|}{|F_n|} \mbox{, for any }s,y\in G\]
\end{lem}

\begin{pr} 
First \[\sum _{x\in yL_n} \phi _n^2(x) = \sum _{x\in yL_n}\frac{|K_n\cap F_nx|}{|F_n|}=\sum _{l\in L_n}\frac{|K_nl\cap F_ny|}{|F_n|}\] after letting $x=yl^{-1}$ with $l\in L_n^{-1}=L_n$. Recall that $G$ has a partition of the form $\cup _{l\in L_n}K_nl$, hence the above sum becomes \[\frac{|(\cup _{l\in L_n}K_nl)\cap F_ny|}{|F_n|} = \frac{|G\cap F_ny|}{|F_n|} = 1\]
For the second part:
\[ \sum_{x\in yL_n} |\phi _n^2(x)-\phi _n^2(sx)| = \sum_{x\in yL_n}\left|\frac{|K_n\cap F_nx|-|K_n\cap F_nsx|}{|F_n|}\right| \] The trick is to cancel out their intersection $|K_n\cap (F_nx\cap F_nsx)|$, thus getting \[ \sum_{x\in yL_n}\left|\frac{|K_n\cap (F_nx\setminus F_nsx)|-|K_n\cap (F_nsx\setminus F_nx)|}{|F_n|}\right| \] This will be less than or equal to \[\sum_{x\in yL_n}\left(\frac{|K_n\cap (F_n\setminus F_ns)x|+|K_n\cap (F_ns\setminus F_n)x|}{|F_n|}\right)\] and by going through the partition argument once more, we end up with \[ \frac{|(F_n\setminus F_ns)y|+|(F_ns\setminus F_n)y|}{|F_n|} = \frac{|F_n\setminus F_ns|+|F_ns\setminus F_n|}{|F_n|} = \frac{|F_n \triangle F_ns|}{|F_n|}\]\qed
\end{pr}

Set $\xi _{ yL_n}=\sum _{x\in yL_n}\phi _n(x)\delta _x$ for every $y\in K_n$. This gives an orthonormal family of $|K_n|$-many vectors in $\ell^2(G)$. Let $P_{ yL_n}$ be the rank-one projection onto $\xi _{ yL_n}$, and $P_n$ be their sum $\sum _{y\in K_n}P_{ yL_n}$. By letting $n\to\infty$, the projection $P_n$ converges to the identity~$I$ on $\ell^2(G)$, as the F\o lner sets $F_n\subset K_n$ exhaust the group $G$. Moreover,

\begin{lem}
$\|[\lambda (s),P_n]\|\to 0$ for all $s\in G$ as $n\to\infty$.
\end{lem}

\begin{pr}
We first observe that \[\|[\lambda (s),P_n]\|\le \|(I-P_n)\lambda (s)P_n\|+\|(I-P_n)\lambda (s^{-1})P_n\|\]for all $s\in G$. To verify, write $I=P_n+(I-P_n)$, then apply the triangle inequality and an involution. We focus on the first of these terms.

Given a unit vector $\xi\in\ell^2(G)$, we have $P_n\xi= \sum _{y\in K_n}\langle\xi ,\xi _{ yL_n}\rangle\xi _{ yL_n}$ with \\$\sum _{y\in K_n}|\langle\xi ,\xi _{ yL_n}\rangle |^2\le 1$. Since \[\lambda (s)\xi _{ yL_n} =\sum _{x\in yL_n}\phi _n(x)\delta _{sx}=\xi _{syL_n}+\sum _{x\in yL_n}(\phi _n(x)-\phi _n(sx))\delta _{sx}\] we obtain: \[\|(I-P_n)\lambda (s)P_n\xi\|^2\le\sum _{x\in yL_n}|\phi _n(x)-\phi _n(sx)|^2\le\sum _{x\in yL_n}|\phi^2_n(x)-\phi^2_n(sx)|\] Now use the second part of the previous lemma (and repeat the argument for the other term) to conclude that \[\|[\lambda (s),P_n]\|\le 2 \sqrt{\frac{|F_n \triangle F_ns|}{|F_n|}}\]\qed
\end{pr}

\begin{rem}
A group $G$ is \emph{maximally almost periodic} if it has a separating family of finite dimensional unitary representations. Every discrete group is the increasing union of its finitely generated subgroups, and it follows from Mal'tsev's Theorem (which states that finitely generated linear groups are residually finite) that a discrete maximally almost periodic group is the increasing union of its residually finite subgroups. Taking into account that quasidiagonality is preserved under increasing unions of C*-algebras, we recover the following theorem of M. Bekka (\cite{1}):
\end{rem}

\begin{thm}
For $G$: discrete countable amenable and maximally almost periodic, the group C*-algebra $C^*(G)$ is quasidiagonal.
\end{thm}

Bekka's theorem is slightly more general: it says that the group C*-algebra is residually finite dimensional (in other words, it admits a block-diagonal representation). 

\section{The main result}

M. Pimsner and D. Voiculescu, in their work on the short exact sequences for K-groups and Ext-groups (\cite{4}), obtained the following result:

\begin{thm}[Pimsner--Voiculescu]
Let $A$ be a unital separable C*-algebra and let $\alpha :\mathbb{Z} \to Aut A$ be a homomorphism, such that there exists a sequence of integers $0\le j_1< j_2<\dotsm$ so that $\lim _{n\to\infty} \|\alpha (j_n)a - a\|=0$ for all $a\in A$. Assume, moreover, that $A$ is quasidiagonal. Then $A\rtimes _\alpha\mathbb{Z}$ is also quasidiagonal.
\end{thm}

Motivated by their method, we prove:

\begin{thm}
Let $G$ be a discrete countable amenable and residually finite group with a sequence of F\o lner sets $F_n$ and tilings of the form $G=K_nL_n$ with $F_n\subset K_n$ for all $n\ge 1$. Let $A$ be a unital separable C*-algebra and let $\alpha :G \to Aut A$ be a homomorphism such that \[\lim _{n\to\infty} \left(\max_{l\in L_n \cap K_nK_n^{-1}F_n}\|\alpha (l)a - a\|\right) = 0\]for all $a\in A$. Assume, moreover, that $A$ is quasidiagonal. Then $A\rtimes _\alpha G$ is also quasidiagonal.
\end{thm}

Let us note that for $G=\mathbb{Z}$, Theorem 7 coincides with Pimsner--Voiculescu's result. For $G=\mathbb{Z}^m$, we have $m$-commuting actions and they should all be almost periodic in the same sense as in the one dimensional case.

\begin{pr}
We adapt the proof of \cite{4} to our more general context.
Let $\rho :A\to \mathcal{B(H)}$ be a faithful quasidiagonal *-representation. Consider the (faithful) regular representation $(\sigma ,I_{\mathcal{H}}\otimes \lambda, \mathcal{H}\otimes\ell^2(G))$, where \[\sigma (a) (h\otimes \delta _x) = \rho (\alpha (x^{-1})a)h\otimes \delta _x\] for $a\in A$. Let $(a_i)^{\infty}_{i=1}$ be a dense subset of all self-adjoint elements $A^{sa}\subset A$. By dropping to a subsequence of tilings $K_nL_n$ of $G$, if necessary, we may assume that $\|\alpha (l)a_i-a_i\|\le\frac{1}{n}$ for all $i\le n$ and all $l\in L_n\cap K_nK_n^{-1}F_n$. Consider $\xi_{ yL_n}$, $P_{ yL_n}$ and $P_n$ as before. Let $Q_n\in \mathcal{B(H)}$ be finite rank projections such that $Q_n\to I_{\mathcal{H}}$ and $\|[\rho (\alpha (y^{-1})a_i),Q_n]\|\le\frac{1}{n}$ for all $i\le n$ and $y\in K_n$. The projections $Q_n\otimes P_n$ converge to the identity~$I$ on $\mathcal{H}\otimes\ell^2(G)$ as $n\to\infty$, and, due to Lemma~4, they commute in the limit with $I_{\mathcal{H}}\otimes \lambda (s)$.

To prove the remaining part, note that for cosets $yL_n\ne zL_n$, we get \[[\sigma (a_i),Q_n\otimes P_{yL_n}]^*[\sigma (a_i),Q_n\otimes P_{zL_n}]=0\]As a result, $[\sigma (a_i),Q_n\otimes P_n]$ is block-diagonal and its norm equals the maximum of the norms of its blocks. Moreover, \[\|[\sigma (a_i),Q_n\otimes P_{yL_n}]\| \le 2\|(I-Q_n\otimes P_{yL_n})\sigma (a_i)(Q_n\otimes P_{yL_n})\|\]Take a unit vector $h\otimes \xi _{yL_n}$ in $Q_n\mathcal{H}\otimes P_n\ell^2(G)$ and compute
\[\sigma (a_i)(h\otimes \xi _{ yL_n}) = \sigma (a_i)(h\otimes \sum _{x\in yL_n}\phi_n(x)\delta_x)=\sum _{x\in yL_n}\phi_n(x)[\rho (\alpha (x^{-1})a_i)h\otimes \delta _x]\] Now write the last quantity as \[\sum _{x\in yL_n}\phi_n(x)[\rho (\alpha (y^{-1})a_i)h\otimes \delta _x]+\sum _{x\in yL_n}\phi_n(x)\left[[\rho (\alpha (x^{-1})a_i)-\rho (\alpha (y^{-1})a_i)]h\otimes \delta _x\right]\] \[ = \rho (\alpha (y^{-1})a_i)h\otimes \xi_{ yL_n}+\sum _{x\in yL_n}\phi_n(x)[\rho (\alpha (x^{-1})a_i-\alpha (y^{-1})a_i)h\otimes \delta _x]\]We first show that the second term has small norm. Indeed, \[\sqrt{\sum _{x\in yL_n}\phi^2_n(x)\|\rho (\alpha (x^{-1})a_i-\alpha (y^{-1})a_i)\|^2} \le \sqrt{\sum _{x\in yL_n}\phi^2_n(x)\|\alpha (x^{-1})a_i-\alpha (y^{-1})a_i\|^2}\] and since $L_n$ is a normal subgroup of $G$, we can write $x\in yL_n=L_n^{-1}y$ with $\phi_n(l^{-1}K_n)=0$ for $l\notin K_nK_n^{-1}F_n$, which gives \[\sqrt{\sum _{l\in L_n\cap K_nK_n^{-1}F_n}\phi^2_n(l^{-1}y)\|\alpha (l)a_i-a_i\|^2}\le\max _{l\in L_n\cap K_nK_n^{-1}F_n}\|\alpha (l)a_i-a_i\|\le\frac{1}{n}\]On the other hand, $(I-Q_n\otimes P_{ yL_n})\rho (\alpha (y^{-1})a_i)h\otimes \xi_{ yL_n} = [\rho (\alpha (y^{-1})a_i),Q_n]h\otimes\xi _{yL_n}$. Therefore, the norm \[\|(I-Q_n\otimes P_{ yL_n})\sigma (a_i)Q_n\otimes P_{ yL_n}\|\le\|[\rho (\alpha (y^{-1})a_i),Q_n]\| + \frac{1}{n} \le \frac{2}{n}\]
Overall $\|[\sigma (a_i),Q_n\otimes P_n]\|\le\frac{4}{n}$ for $i\le n$. Hence, $\lim _{n\to\infty}\|[\sigma (a_i),Q_n\otimes P_n]\|=0$\qed
\end{pr}

Recall that every finitely generated subgroup of a discrete maximally almost periodic group is residually finite.

\begin{cor}
Assume $G$ is a discrete countable amenable and maximally almost periodic group. Let $\alpha$ be its action on a unital separable quasidiagonal C*-algebra $A$. If the restriction of $\alpha$ to every finitely generated subgroup of $G$ is almost periodic, in the sense of Theorem~7, then the crossed product $A\rtimes_\alpha G$ is quasidiagonal.
\end{cor}

An application of Theorem~7 appears in \cite{3}. There, we consider a broad family of crossed products of the form $C(\tilde{G})\rtimes G$, where $\tilde{G}$ denotes a profinite completion of $G$. We call such crossed products the \emph{generalized Bunce--Deddens algebras}, in view of the well-known fact that the classical Bunce--Deddens algebras can be written in the form $C(\tilde{\mathbb{Z}})\rtimes \mathbb{Z}$. We then prove that the generalized Bunce--Deddens algebras are quasidiagonal, and also that they are unital separable simple nuclear C*-algebras, of real rank zero, stable rank one, with comparability of projections and a unique trace.

\begin{akn}
The results presented here are part of the author's doctoral thesis at Purdue University. I~am indebted to Marius Dadarlat for tons of advice and support and I~would like to thank the Fields Institute and George Elliott for their hospitality during Fall 2007; also Eberhard Kirchberg and Chris Phillips for many helpful discussions. Last but not least, many thanks to Larry Brown for comments on an earlier draft of this paper.
\end{akn}


\end{document}